\documentclass[review]{elsarticle}
\usepackage{amsfonts}
\usepackage{hyperref}
\usepackage{amssymb}
\usepackage{amsmath}
\usepackage{graphicx}
\usepackage{float}
\usepackage[export]{adjustbox}
\usepackage{subfigure}
\usepackage{geometry}
\usepackage{epstopdf}
\usepackage{caption}
\usepackage{booktabs}
\usepackage{bigstrut,multirow,rotating}
\usepackage{xcolor}  
\usepackage{tikz}  
\usepackage{diagbox}
\usepackage{multirow}
\usepackage{array}
\usepackage{longtable}
\usepackage{graphicx}
\usepackage{algpseudocode}
\usepackage{algorithmicx,algorithm}
\usepackage{float}

\usepackage[export]{adjustbox}
\usepackage{subfigure}
\usepackage{geometry}
\usepackage{epstopdf}
\usepackage{caption}
\usepackage{geometry}
\usepackage{amsthm}

\biboptions{sort&compress}
\usetikzlibrary{arrows,shapes,chains}
\begin{document}

\begin{frontmatter}

\title{An  Image Noise Level Estimation Based on Tensor T-Product}

\author[mymainaddress]{Hanxin Liu}

\author[mymainaddress]{Yisheng Song\corref{mycorrespondingauthor}}

\cortext[mycorrespondingauthor]{Corresponding author}
\ead{yisheng.song@cqnu.edu.cn}

\address[mymainaddress]{School of Mathematical Sciences, Chongqing Normal University, Chongqing, 401331.\\ Email: \em 2531417503@qq.com (Liu); yisheng.song@cqnu.edu.cn (Song)}

\begin{abstract}
Currently, the noise level of color images is estimated by many algorithms through separate selection of each page of the third-order tensor using sliding blocks of size ${M_1} \times {M_1}$. The data structure of the tensor is disrupted by this method, leading to errors in the estimation results. In order not to disrupt the data structure of the tensor, we directly select the tensor using a sliding block of size ${M_1} \times {M_1} \times 3$ and then re-arrange it. The newly obtained tensor is decomposed into a block diagonal matrix form through T-product. It is demonstrated that the eigenvalues of this matrix are related to the noise level of the color image. Then train the relationship coefficients through learning methods, thereby obtaining the estimated noise level. The effectiveness of the algorithm was verified through numerical experiments, and it also achieved high estimation accuracy.

\end{abstract}


\end{frontmatter}


\section{Introduction}

The image noise level indicates the degree to which an image is contaminated by noise, and serving as a crucial parameter in the image processing process. For instance, in processes such as blind image denoising\cite{ref1}, blind image restoration\cite{ref2}, and blind image deblurring\cite{ref3}, estimating the noise level beforehand is necessary. Hence, it is important to devising an accurate noise level estimation method.

Over the past decade, the research on noise level estimation has been a hot topic in image processing. There are many methods for estimating noise levels, such as filtering-based\cite{ref4}, transformation-based\cite{ref5} and block-based\cite{ref6} methods. This study mainly focuses on the block-based method for estimating noise levels. This method divides the image into overlapping small blocks through sliding windows, selects the blocks with consistent statistical data (such as variance and kurtosis), forming uniform blocks or weakly textured blocks, and finally calculates the noise level of the image in these selected blocks. 

A block-based noise level estimation method was first proposed by Pyatykh et al.\cite{ref7}. They estimated the noise level using image block principal component analysis. Even when the noise image contains very few uniform areas, this typical method can still estimate the noise level of the image with a certain degree of accuracy. Based on this method, Liu et al.\cite{ref8}  proposed an image noise level estimation algorithm based on weak texture blocks. They select weakly textured image blocks from images with rich texture information to estimate the noise level. This algorithm can effectively estimate the noise level of images with rich texture information. The above two algorithms both calculate the covariance matrix of the selected image block, and use the minimum eigenvalue of the covariance matrix as the noise level of the image. The color image is regarded as a third-order tensor. Each slice is treated as a matrix, and it is divided into three matrices to perform the above operations respectively. It was pointed out by Heng et al.\cite{ref9} that the noise level estimated by the above method is a univariate function of the minimum eigenvalue, which often leads to underestimation. Subsequently, based on the calculation of the eigenvalues of the covariance matrix, they proposed to compare the median and the mean of the eigenvalues. If the median and the mean are not equal, then remove the largest eigenvalue, and calculate the median and mean of the remaining eigenvalues. Repeat this process until the median and the mean are equal, at which point the final estimation result is obtained. 

Based on the above approach of using a single eigenvalue as the noise level, some scholars have taken into account other factors. Fang et al.\cite{ref10} conducted a further analysis of the eigenvalues in 2019. They believe that using the minimum eigenvalue as the noise variance would lead to underestimation, while using the mean of the eigenvalues would result in overestimation. They perform a linear fitting on the results of an overestimation and an underestimation, and by proving that the estimated noise level has the following linear relationship ${\sigma ^2} = ({d_1}\sigma _f^2 + {d_2}\sigma _w^2)/({d_1} + {d_2})$, where $\sigma _f^2$ is an underestimated result, $\sigma _w^2$ is an overestimated result. In 2020, Jiang et al.\cite{ref11} took into account the factors of the number of image blocks and the size of the sliding blocks. Their further research reveals that the relationship between noise level and the minimum eigenvalue ${\lambda _{\min }}$, the number of image blocks $s$, and the size of sliding blocks $w$ is as follows ${\sigma ^2} = {\lambda _1}/(1 - 1.8606\sqrt {({w^2} - 2)/s} )$. Liu et al.\cite{ref12} estimated the noise level by fitting multiple feature values on the training set through a learning method.

When estimating the noise level of color images using the above method, the sliding blocks are used to process the three matrices separately, and three covariance matrices are constructed. The noise level is estimated by calculating the eigenvalues of these matrices. Color images can be regarded as a third-order tensor, grayscale image is a second-order tensor, that is, a matrix. Splitting a third-order tensor into three matrices will disrupt the structure of the tensor and may lead to errors in the estimation results. 

From the perspective of tensors, directly calculating the eigenvalues of higher-order tensors is NP-hard. Unlike matrices, the eigenvalues of tensors have many definitions, including H-eigenvalues\cite{ref13}, Z-eigenvalues\cite{ref14}, M-eigenvalues\cite{ref15}, D-eigenvalues\cite{ref16}, and B-eigenvalues\cite{ref17}, etc. Since Professor Qi Liqun proposed the definition of tensor special H-eigenvalues in 2005\cite{ref13}, the calculation of tensor eigenvalues has been a hot topic in this field. As proposed by Qi et al.\cite{ref18}, the Z-eigenvalues of a third-order three-dimensional tensor can be directly calculated using the orthogonal transformation method. Kofidis et al. proposed the higher-order power method (HOPM)\cite{ref19} and the symmetric higher-order power method (S-HOPM)\cite{ref19} for calculating the eigenvalues of large-scale tensors. Subsequently, Kolda et al. added displacement parameters to the S-HOPM algorithm and controlled the convergence of the algorithm through these displacement parameters. This algorithm is called the Shifted Symmetric High-Order Power Method (SS-HOPM)\cite{ref20}. This algorithm is overly dependent on the selection of parameters. Subsequently, Kolda et al. proposed an adaptive displacement method to address the issue of excessive reliance on parameters. This method adaptively selects the displacement parameters based on the positive and negative definiteness of the Hessian matrix of the objective function, and is called the  Generalized Eigenproblem Adaptive Power (GEAP)\cite{ref21}.

The above algorithms all calculate the eigenvalues of some specific tensors. It is impossible to directly calculate the eigenvalues of the third-order tensor corresponding to the color image. Therefore, a decomposition method such as CP decomposition\cite{ref22} and Tucker decomposition\cite{ref23} is needed to decompose the tensor into a form similar to a matrix. Kilmer et al.\cite{ref24} are the first to propose the concepts of T-SVD decomposition and T-product. Under the definition of T-product, this decomposition decomposes the higher-order tensor into a form similar to the outer product of matrices. Based on the idea of T-SVD, in this paper, under the definition of T-product, the third-order tensor is transformed into a matrix form through calculation. Calculate the eigenvalues of the covariance matrix of the new matrix and thereby estimate the noise level of the color image. The main innovations of this paper are as follows: 
\begin{itemize}
\item A new noise level estimation model based on tensor decomposition is proposed. The third-order tensor is decomposed into the form of a covariance matrix according to the definition of the T-product, and then the noise level of the color image is estimated, and good experimental results were achieved in the experiment.
\item Through theoretical analysis, it has been proved that the multiple eigenvalues of the matrix obtained by T-product of the third-order tensor have a direct relationship with the noise level of the color image.

\end{itemize}

The rest of the paper is structured as follows: Basic knowledge, including the traditional noise level estimation method and T - product, is introduced in Section \ref{2}. The proposed noise level estimation algorithm is described in Section \ref{3}. The results of this paper, which are compared with those of other advanced algorithms, are presented in Section \ref{4}. The algorithm in this paper is summarized and the future work is briefly described as well in Section \ref{5}.



\section{Background knowledge}\label{2}
\subsection{Image noise level estimation}

For an observed image $y$ that contains additive Gaussian white noise, its model can be expressed in the following form
\begin{equation}\label{M1}
{y} = {x} + {e},
\end{equation}
where $x$  is the noiseless image patch,  $e$ is the signal-independent additive white Gaussian noise with mean value of 0 and variance of ${\sigma ^2}$. Assume that the size of the observed image $y$ is ${S_1} \times {S_2} \times c$. Divide $y$ into $s = ({S_1} - {M_1} + 1) \times ({S_2} - {M_1} + 1) \times c$ patches using ${M_1} \times {M_1}$ -sized sliding window, then reorder each patche into a column vector of size $M_1^2 \times 1\ $.  So that an observation image $y$  can be represented as ${Y_s} = \{ {y_t}\} _{t = 1}^s \in $$\mathbb{R}^{M_1^2 \times s}$ by $s$ image patches. The covariance matrix ${\Sigma _y}$ of the noise image $y$ is defined as:
\begin{equation}\label{M2}
{\Sigma _y} = \frac{1}{s}\sum\limits_{i = 1}^s {({y_i} - u){{({y_i} - u)}^T}} ,
\end{equation}
where column vector $u$ is the average value of the data set $\{ {y_i}\} $. Covariance matrix ${\Sigma _y}$ satisfies the following assumption:

\textbf{Assumption 1:} Covariance matrix ${\Sigma _y}$ follows the gamma distribution with shape parameter $(s - 1)/2$ and scale parameter $(s - 1)/2{\sigma ^2}$ :  
 \begin{equation}\label{M3}
{\Sigma _y} \sim \gamma ((s - 1)/2,(s - 1)/2{\sigma ^2}),
\end{equation}
where $\gamma $ is the gamma distribution, expectation is ${\sigma ^2}$, variance is $2{\sigma ^4}/(s - 1)$. 

Under the \textbf{Assumption 1}, the minimum eigenvalue of the observed image covariance matrix and the original image covariance matrix satisfies the following relationship
\begin{equation}\label{M4}
{\lambda _{\min }}({\Sigma _y}) = {\lambda _{\min }}({\Sigma _x}) + {\sigma ^2},
\end{equation}
where ${\Sigma _y}$ is the covariance matrix of noisy image patch ${y_i}$, ${\Sigma _x}$ is the covariance matrix of noiseless image patch ${x_i}$, and ${\lambda _{\min }}(\Sigma )$ represents the minimum eigenvalue of matrix $\Sigma $. The noise level ${\sigma ^2}$ can be calculated, if the minimum eigenvalue of  ${\Sigma _x}$ is known. As we know, the minimum eigenvalue of the weak texture patches' covariance matrix is 0. Therefore, we select the image blocks with weak texture from the images to estimate the noise level. Liu et al.\cite{ref8} used the iterative threshold $\tau $ to obtain weak texture patches, and the threshold $\tau $ is defined as follows:
\begin{equation}\label{M5}
\tau  = {\sigma ^2}{F^{ - 1}}(\delta ,\frac{N}{2},\frac{2}{N}tr(D_h^T{D_h} + D_v^T{D_v})),
\end{equation}
where ${F^{ - 1}}(\delta ,\alpha ,\beta )$ is gamma cumulative distribution inverse function with shape parameter $\alpha $ and scale parameter $\beta$, ${D_h}$ and ${D_v}$ represent horizontal and vertical operator matrix respectively. When the maximum eigenvalue of an image patch covariance matrix is less than this threshold $\tau $, the image patch is the weak texture patch required for noise level estimation. According to the above method, the noise level of the observed image $y$ can be estimated as:
\begin{equation}\label{M6}
{\hat \sigma ^2} = g({\lambda _{\min }}).
\end{equation}
\subsection{Tensor T-product}

If $v = {[\begin{array}{*{20}{c}}
{{v_0}}&{{v_1}}&{{v_2}}&{{v_3}}
\end{array}]^T}$, then 
\[circ(v) = \left[ {\begin{array}{*{20}{c}}
{{v_0}}&{{v_3}}&{{v_2}}&{{v_1}}\\
{{v_1}}&{{v_0}}&{{v_3}}&{{v_2}}\\
{{v_2}}&{{v_1}}&{{v_0}}&{{v_3}}\\
{{v_3}}&{{v_2}}&{{v_1}}&{{v_0}}
\end{array}} \right],\]
is a circulant matrix. A color image can be regarded as a third-order tensor ${\cal A}$ of size ${n_1} \times {n_2} \times 3$. Just like the creation of a circulant matrix, it is possible to create a block circulant matrix from the slices of a tensor. Then a tensor ${\cal A} \in {\mathbb{R}^{{n_1} \times {n_2} \times {n_3}}}$ can be transformed into the following form
\[circ({\cal A}) = \left[ {\begin{array}{*{20}{c}}
{{\cal A}{_1}}&{{\cal A}{_{{n_3}}}}& \cdots &{{\cal A}{_2}}\\
{{\cal A}{_2}}&{{\cal A}{_1}}& \cdots &{{\cal A}{_3}}\\
 \vdots & \ddots & \ddots & \vdots \\
{{\cal A}{_{{n_3}}}}&{{\cal A}{_{{n_3} - 1}}}& \cdots &{{\cal A}{_1}}
\end{array}} \right],\]
where ${{\cal A}_i} = {\cal A}(:,:,i)$ for $i = 1,2, \cdots ,{n_3}$.

\textbf{Definition:} Let a tensor ${\cal A} \in {\mathbb{R}^{{n_1} \times {n_2} \times {n_3}}}$ and a tensor ${\cal B} \in {\mathbb{R}^{{n_2} \times {n_4} \times {n_3}}}$. Then the t-product ${\cal A} * {\cal B} \in {\mathbb{R}^{{n_1} \times {n_4} \times {n_3}}}$  
\begin{equation}\label{M7}
{\cal A}*{\cal B} = fold\left( {circ({\cal A}) \cdot MatVec({\cal B})} \right).
\end{equation}
where
\[MatVec({\cal B}) = \left[ {\begin{array}{*{20}{c}}
{{{\cal B}_1}}\\
{{{\cal B}_2}}\\
 \vdots \\
{{{\cal B}_{{n_3}}}}
\end{array}} \right],\]
the operation that takes $MatVec(\cal A)$ back to tensor form is the fold command: $fold(MatVec(\cal A)) = \cal A$.

\section{Proposed algorithm}\label{3}
\subsection{Model in this paper}
For a color image, the traditional method of calculating the covariance matrix will disrupt the data structure of the tensor. This paper operates on the entire tensor, using a sliding window of size ${M_1} \times {M_1} \times 3$ to extract $s = ({n_1} - {M_1} + 1) \times ({n_2} - {M_1} + 1) $ block tensors from a third-order tensor of size ${n_1} \times {n_2} \times 3$. Rearrange each slice of the block tensor into a column vector, so that each block tensor is reorganized into a matrix of size ${M_1}^2 \times 3$. Combine each matrix obtained from the block tensors to form a third-order tensor ${\cal A}$ of size ${M_1}^2 \times s \times 3$. $\{ y_i^j\} $ represent the column vector of the j-th slices in the i-th column of ${\cal A}$. ${u^j}$ represent the mean value of $\{ y_i^j\} $. Then, the covariance matrix of each slice of tensor ${\cal A}$ can be expressed as
\begin{equation}\label{M8}
{\Sigma _{\cal A}{^j}} = \frac{1}{s}\sum\limits_{i = 1}^s {({y_i^j} - u){{({y_i^j} - u)}^T}}.
\end{equation}
Reorganize the covariance matrix of each slice into a third-order tensor ${\cal B} = ({\Sigma _{{{\cal A}^j}}})$(j=1,2,3). 

In order to decompose the aforementioned tensor into a form similar to a covariance matrix, we introduce a tensor ${\cal I}_i$, where the i-th slice is a unit matrix $E$. Then the third-order tensor ${\cal B} \in {\mathbb{R}^{{M_1}^2 \times {M_1}^2 \times {3}}}$ is transformed as follows to obtain the matrix B. 

\[B = unfold({\cal B}*{{\cal I}_2}) = \left( {\begin{array}{*{20}{c}}
{\left[ {\begin{array}{*{20}{c}}
{\Sigma _{{{\cal A}^1}}}&{\Sigma _{{{\cal A}^3}}}&{\Sigma _{{{\cal A}^2}}}\\
{\Sigma _{{{\cal A}^2}}}&{\Sigma _{{{\cal A}^1}}}&{\Sigma _{{{\cal A}^3}}}\\
{\Sigma _{{{\cal A}^3}}}&{\Sigma _{{{\cal A}^2}}}&{\Sigma _{{{\cal A}^1}}}
\end{array}} \right]}&{\left[ {\begin{array}{*{20}{c}}
0\\
{{E_2}}\\
0
\end{array}} \right]}
\end{array}} \right) = \left[ {\begin{array}{*{20}{c}}
{\Sigma _{{{\cal A}^3}}}\\
{\Sigma _{{{\cal A}^1}}}\\
{\Sigma _{{{\cal A}^2}}}
\end{array}} \right].\]
We use the operator $bdiag( \cdot )$ to represent the transformation of matrix B into block diagonal matrix $bdiag(B)$.
\[bdiag(B) = \left[ {\begin{array}{*{20}{c}}
{{\Sigma _{{{\cal A}^3}}}}&{}&{}\\
{}&{{\Sigma _{{{\cal A}^1}}}}&{}\\
{}&{}&{{\Sigma _{{\cal A}^2}}}
\end{array}} \right].\]

Since ${{\Sigma _{{{\cal A}^j}}}}$ is the covariance matrix composed of each slice of ${\cal A}$, and ${{\Sigma _{{{\cal A}^j}}}}$ satisfies \textbf{Assumption 1}, then the following theorem holds.

\textbf{Theorem 2:} If ${\lambda _{1}},{\lambda _{2}}, \cdot  \cdot  \cdot ,{\lambda _{r}}$ are $r$ eigenvalues of the matrix $bdiag(B)$ from small to large. Under \textbf{Assumption 1}, there exists ${\theta _1},{\theta _2}, \cdot  \cdot  \cdot ,{\theta _n}(n \leqslant r)$, satisfying ${\theta _1} + {\theta _2} +  \cdot  \cdot  \cdot  + {\theta _n} = 1$, such that: 
 \begin{equation}\label{M9}
{\sigma ^2} - \sqrt {2r} {\sigma ^2}/\sqrt {(s - 1)}  \le {\theta _1}{\lambda _1} + {\theta _2}{\lambda _2} +  \cdots  + {\theta _n}{\lambda _n} \le {\sigma ^2} + \sqrt {2r} {\sigma ^2}/\sqrt {(s - 1)} 
\end{equation}
\begin{proof}
We know that $\{ {\lambda _{bdiag(B)}}\}  = \{ {\lambda _{{\Sigma _{{{\cal A}^3}}}}}\}  \cup \{ {\lambda _{{\Sigma _{{{\cal A}^2}}}}}\}  \cup \{ {\lambda _{{\Sigma _{{{\cal A}^1}}}}}\} $. When we are discussing the eigenvalues of $bdiag(B)$, we can successively discuss the eigenvalues of ${{\Sigma _{{{\cal A}^3}}}}$, ${{\Sigma _{{{\cal A}^1}}}}$ and ${{\Sigma _{{{\cal A}^2}}}}$.

First, we discuss the eigenvalues of ${{\Sigma _{{{\cal A}^3}}}}$. Suppose $ \lambda _1^3,\lambda _2^3, \cdots \lambda _{{r_3}}^3$ are the eigenvalues of ${{\Sigma _{{{\cal A}^3}}}}$, any real number ${\beta _j}(j = 1,2, \cdot  \cdot  \cdot ,r_3)$, 
\[\begin{array}{c}
\begin{aligned} 
\sum\limits_{j = 1}^{{r_3}} {\bar \beta (\bar \lambda^3  - {\lambda _j^3})}&= {r_3}\bar \beta \bar \lambda^3  - \bar \beta ({\lambda _1^3} + {\lambda _2^3} +  \cdots  + {\lambda _{{r_3}}^3})\\
 & = {r_3}\bar \beta \frac{{({\lambda _1^3} + {\lambda _2^3} +  \cdots  + {\lambda _{{r_3}}^3})}}{{{r_3}}} - \bar \beta ({\lambda _1^3} + {\lambda _2^3} +  \cdots  + {\lambda _{{r_3}}^3})\\
 &= 0
\end{aligned} 
\end{array},\]
where $\bar \beta  = \tfrac{1}{r_3}\sum {{\beta _j}}$, according to Cauchy-Schwartz inequality, we have 
 \begin{equation}\label{M10}
\begin{array}{c}
\begin{aligned} 
\left| {\sum\limits_{j = 1}^{{r_3}} {{\beta _j}(\bar \lambda^3  - {\lambda _j^3})} } \right| &= \left| {\sum\limits_{j = 1}^{{r_3}} {{\beta _j}(\bar \lambda^3  - {\lambda _j^3})}  - \sum\limits_{j = 1}^{{r_3}} {\bar \beta (\bar \lambda^3  - {\lambda _j^3})} } \right|\\
 &= \left| {\sum\limits_{j = 1}^{{r_3}} {({\beta _j} - \bar \beta )(\bar \lambda^3  - {\lambda _j^3})} } \right|\\
 & \le {\left| {\sum\limits_{j = 1}^{{r_3}} {{{({\beta _j} - \bar \beta )}^2}} \sum\limits_{j = 1}^{{r_3}} {{{(\bar \lambda^3  - {\lambda _j^3})}^2}} } \right|^{\frac{1}{2}}},
\end{aligned} 
\end{array}
\end{equation}
where $\bar \beta  = \tfrac{1}{r}\sum {{\beta _j}}$, according to \textbf{Assumption 1},
\begin{equation}\label{M11}
E({\Sigma _{{{\cal A}^3}}}) = \bar \lambda^3  = {\sigma ^2},
\end{equation}
\begin{equation}\label{M12}
D({\Sigma _{{{\cal A}^3}}}) = 2r{{\sigma ^4}/(s - 1)},
\end{equation}
we obtain
\begin{equation}\label{M13}
\begin{array}{c}
\begin{aligned} 
\sum {({{\bar \lambda }^3} - \lambda _j^3)^2}  
&\leqslant \sum {(\lambda _j^3)^2} - r_3{({{\bar \lambda }^3})^2}\\ 
&\leqslant \sum (E{({\Sigma _{{{\cal A}^3}}})^2}) - r_3{({{\bar \lambda }^3})^2}\\ 
&= \sum (D({\Sigma _{{{\cal A}^3}}})) + \sum {(E({\Sigma _{{{\cal A}^3}}}))^2} - r_3{({{\bar \lambda }^3})^2}\\
&= 2r_3{\sigma ^4}/(s - 1) + r_3{\sigma ^2} - r_3{\sigma ^2}\\
&= 2r_3{\sigma ^4}/(s - 1)
\end{aligned} 
\end{array}.
\end{equation}
Moreover,

\begin{equation}\label{M14}
\left| {\sum {{\beta _j}(\bar \lambda^3  - {\lambda _j^3})} } \right| \leqslant {\left| 2r_3{{\sigma ^4}/(s - 1)}{\sum {{{({\beta _{j}} - \bar \beta )}^2}} } \right|^{\frac{1}{2}}},
\end{equation}
 let ${\beta _1} = 1$ and ${\beta _j} = 0(j \ne 1)$, combine Eq. (\ref{M11}) and (\ref{M14}) the above inequality can be rewritten as
\begin{equation}\label{M15}
\left| {{\sigma ^2} - \lambda _1^3} \right| \le \sqrt {2r} {\sigma ^2}/\sqrt {(s - 1)}.
\end{equation}
Therefore,
\begin{equation}\label{M16}
{\sigma ^2} - \sqrt {2r_3} {\sigma ^2}/\sqrt {(s - 1)}  \le \lambda _1^3 \le {\sigma ^2} + \sqrt {2r_3} {\sigma ^2}/\sqrt {(s - 1)}.
\end{equation}
In a similar way, let ${\beta _2} = 1$ and ${\beta _j} = 0(j \ne 2); \cdot  \cdot  \cdot ; {\beta _n} = 1$ and ${\beta _j} = 0(j \ne n)$, 
\begin{equation}\label{M17}
\begin{gathered}
  {\sigma ^2} - \sqrt {2r_3} {\sigma ^2}/\sqrt {(s - 1)}  \le \lambda _2^3 \le {\sigma ^2} + \sqrt {2r_3} {\sigma ^2}/\sqrt {(s - 1)}  \\ 
  {\sigma ^2} - \sqrt {2r_3} {\sigma ^2}/\sqrt {(s - 1)}  \le \lambda _3^3 \le {\sigma ^2} + \sqrt {2r_3} {\sigma ^2}/\sqrt {(s - 1)} \\ 
   \vdots  \\ 
  {\sigma ^2} - \sqrt {2r_3} {\sigma ^2}/\sqrt {(s - 1)}  \le \lambda _n^3 \le {\sigma ^2} + \sqrt {2r_3} {\sigma ^2}/\sqrt {(s - 1)}\\ 
\end{gathered}.
\end{equation}

Then, we discuss the eigenvalues of ${{\Sigma _{{{\cal A}^1}}}}$ and ${{\Sigma _{{{\cal A}^2}}}}$.The same as the proof of ${{\Sigma _{{{\cal A}^3}}}}$, we have
\begin{equation}\label{M18}
\begin{gathered}
  {\sigma ^2} - \sqrt {2r_1} {\sigma ^2}/\sqrt {(s - 1)}  \le \lambda _1^1 \le {\sigma ^2} + \sqrt {2r_1} {\sigma ^2}/\sqrt {(s - 1)}  \\ 
  {\sigma ^2} - \sqrt {2r_1} {\sigma ^2}/\sqrt {(s - 1)}  \le \lambda _2^1 \le {\sigma ^2} + \sqrt {2r_1} {\sigma ^2}/\sqrt {(s - 1)} \\ 
   \vdots  \\ 
  {\sigma ^2} - \sqrt {2r_1} {\sigma ^2}/\sqrt {(s - 1)}  \le \lambda _n^1 \le {\sigma ^2} + \sqrt {2r_1} {\sigma ^2}/\sqrt {(s - 1)}\\ 
\end{gathered},
\end{equation}
\begin{equation}\label{M19}
\begin{gathered}
  {\sigma ^2} - \sqrt {2r_2} {\sigma ^2}/\sqrt {(s - 1)}  \le \lambda _1^2 \le {\sigma ^2} + \sqrt {2r_2} {\sigma ^2}/\sqrt {(s - 1)}  \\ 
  {\sigma ^2} - \sqrt {2r_2} {\sigma ^2}/\sqrt {(s - 1)}  \le \lambda _2^2 \le {\sigma ^2} + \sqrt {2r_2} {\sigma ^2}/\sqrt {(s - 1)} \\ 
   \vdots  \\ 
  {\sigma ^2} - \sqrt {2r_2} {\sigma ^2}/\sqrt {(s - 1)}  \le \lambda _n^2 \le {\sigma ^2} + \sqrt {2r_2} {\sigma ^2}/\sqrt {(s - 1)}\\ 
\end{gathered},
\end{equation}
where $ \lambda _1^1,\lambda _2^1, \cdots \lambda _{{r_1}}^1$ are the eigenvalues of ${{\Sigma _{{{\cal A}^1}}}}$, $ \lambda _1^2,\lambda _2^2, \cdots \lambda _{{r_2}}^2$ are the eigenvalues of ${{\Sigma _{{{\cal A}^2}}}}$.

Finally, let 
\[\begin{array}{c}
\begin{aligned} 
{\lambda _1} &= \min \{ \lambda _1^3,\lambda _2^3, \cdots ,\lambda _{{r_3}}^3,\lambda _1^1,\lambda _2^1, \cdots ,\lambda _{{r_1}}^1,\lambda _1^2,\lambda _2^2, \cdots ,\lambda _{{r_2}}^2\} \\
{\lambda _2} &= \min \{ \{ \lambda _1^3,\lambda _2^3, \cdots ,\lambda _{{r_3}}^3,\lambda _1^1,\lambda _2^1, \cdots ,\lambda _{{r_1}}^1,\lambda _1^2,\lambda _2^2, \cdots ,\lambda _{{r_2}}^2\} /\{ {\lambda _1}\} \} \\
 \vdots \\
{\lambda _n} &= \min \{ \{ \lambda _1^3,\lambda _2^3, \cdots ,\lambda _{{r_3}}^3,\lambda _1^1,\lambda _2^1, \cdots ,\lambda _{{r_1}}^1,\lambda _1^2,\lambda _2^2, \cdots ,\lambda _{{r_2}}^2\} /\{ {\lambda _1},{\lambda _2}, \cdots ,{\lambda _{n - 1}}\} \} 
\end{aligned} 
\end{array},\]
and $r = \max \{ {r_1},{r_2},{r_3}\}$, based on the above derivation process, we obtain
\begin{equation}\label{M20}
\begin{gathered}
  {\sigma ^2} - \sqrt {2r} {\sigma ^2}/\sqrt {(s - 1)}  \le \lambda _1 \le {\sigma ^2} + \sqrt {2r} {\sigma ^2}/\sqrt {(s - 1)}  \\ 
  {\sigma ^2} - \sqrt {2r} {\sigma ^2}/\sqrt {(s - 1)}  \le \lambda _2 \le {\sigma ^2} + \sqrt {2r} {\sigma ^2}/\sqrt {(s - 1)} \\ 
   \vdots  \\ 
  {\sigma ^2} - \sqrt {2r} {\sigma ^2}/\sqrt {(s - 1)}  \le \lambda _n \le {\sigma ^2} + \sqrt {2r} {\sigma ^2}/\sqrt {(s - 1)}\\ 
\end{gathered},
\end{equation}
So, construct a convex combination, for all ${\theta _1},{\theta _2}, \cdots,{\theta _n}\in(0,1)$ with ${\theta _1} + {\theta _2} +  \cdots  + {\theta _n} = 1$  such that:
 \begin{equation}\label{M22}
{\sigma ^2} - \sqrt {2r} {\sigma ^2}/\sqrt {(s - 1)}  \le {\theta _1}{\lambda _1} + {\theta _2}{\lambda _2} +  \cdots  + {\theta _n}{\lambda _n} \le {\sigma ^2} + \sqrt {2r} {\sigma ^2}/\sqrt {(s - 1)}. 
\end{equation}
The above theorem is proved.
\end{proof}

From \textbf{Theorem 2}, we can see that there exists ${\theta _0} \in [ - \sqrt {2r} {\sigma ^2}/\sqrt {(s - 1)} ,\sqrt {2r} {\sigma ^2}/\sqrt {(s - 1)} ]$ such that 
 \begin{equation}\label{M23}
{{\hat \sigma }^2} = {\theta _0} + {\theta _1}{\lambda _1} + {\theta _2}{\lambda _2} +  \cdots {\theta _n}{\lambda _n},
\end{equation}
where ${{\hat \sigma }^2}$ is the estimated value for the noise level of the observed image. $ \lambda _1,\lambda _2, \cdots \lambda _n$ are the $n$ minimum eigenvalues of matrix $bdiag(B)$.

\subsection{Algorithm flow}

To find the values of ${\theta _j}(j = 0,1, \cdot  \cdot  \cdot ,n)$, learning algorithm is adopted in this paper. $M$ observed images with known noise levels were used as the training set, and the model (\ref{M23}) was input into the training set to calculate the model parameters ${\theta _j}(j = 0,1, \cdot  \cdot  \cdot ,n)$. Specifically speaking, let ${\lambda _0}=1$, for $M$ images, the loss function is constructed as follows: 
\begin{equation}\label{M24}
J({\theta _0},{\theta _1}, \cdot  \cdot  \cdot ,{\theta _n}) = \frac{1}{{2 \times M}}\sum\limits_{i = 1}^M {{{(f(\lambda _0^{(i)},\lambda _{1}^{(i)}, \cdot  \cdot  \cdot ,\lambda _{n}^{(i)}) - {\sigma ^{(i)}})}^2}}, 
\end{equation}
where $\lambda _1^{(i)},\lambda _2^{(i)}, \cdots ,\lambda _n^{(i)}$ are the $n$ smallest eigenvalues of the matrix $bdiag(B)$ for the i-th image, $f({\lambda _0},{\lambda _1}, \cdots ,{\lambda _n}) = {\theta _0}{\lambda _0} + {\theta _1}{\lambda _1} + {\theta _2}{\lambda _2} +  \cdots {\theta _n}{\lambda _n}$. The value of ${\theta _j}(j = 0,1, \cdot  \cdot  \cdot ,n)$ are solved by gradient descent method. The algorithm termination distance $\varepsilon$  and the step size $\alpha$. The gradient expression of the loss function with ${\theta _j}(j = 0,1, \cdot  \cdot  \cdot ,n)$ is:
\begin{equation}\label{M25}
\frac{\partial }{{\partial {\theta _j}}}J({\theta _0},{\theta _1}, \cdot  \cdot  \cdot ,{\theta _n})
 = \frac{1}{M}\sum\limits_{i = 1}^M {(f(\lambda _0^{(i)},\lambda _{1}^{(i)}, \cdot  \cdot  \cdot ,\lambda _{n}^{(i)}) - {\sigma ^{(i)}})\lambda _{ j}^{(i)}} ,
\end{equation}
multiply the gradient of the loss function by the step size $\alpha$ to obtain the descent distance
\begin{equation}\label{M26}
{\alpha}\frac{\partial }{{\partial {\theta _j}}}J({\theta _0},{\theta _1}, \cdot  \cdot  \cdot ,{\theta _n})
 = \frac{\alpha}{M}\sum\limits_{i = 1}^M {(f(\lambda _0^{(i)},\lambda _{1}^{(i)}, \cdot  \cdot  \cdot ,\lambda _{n}^{(i)}) - {\sigma ^{(i)}})\lambda _{ j}^{(i)}} .
\end{equation}
If the descent distance does not meet ${\alpha}\frac{\partial }{{\partial {\theta _j}}}J({\theta _0},{\theta _1}, \cdot  \cdot  \cdot ,{\theta _n}) \le \varepsilon$, then the iterative formula for ${\theta _j}(j = 0,1, \cdot  \cdot  \cdot ,n)$ is
\begin{equation}\label{M27}
\begin{aligned} \theta_{j} &=\theta_{j}-\alpha \frac{\partial}{\partial \theta_{j}} J\left(\theta_{0}, \theta_{1}, \cdots, \theta_{n}\right) \\ &=\theta_{j}-\frac{\alpha}{M} \sum_{i=1}^{M}\left(f\left(\lambda_{0}^{(i)}, \lambda_{1}^{(i)}, \cdots, \lambda_{n}^{(i)}\right)-\sigma^{(i)}\right) \lambda_{j}^{(i)}. \end{aligned}
\end{equation}

The details of the algorithm in this paper are shown in \textbf{Algorithm 1}. For a color image ${\cal A}$ of size $n_1\times n_2 \times 3$, using a sliding window of size ${M_1} \times {M_1} \times 3$ to extract $s = ({n_1} - {M_1} + 1) \times ({n_2} - {M_1} + 1) $ block tensors ${w_{(k)}}$. The block tensors ${w_{(k)}}$ are rearranged into matrices of size ${M_1^2} \times 3$. Reorganize these matrices into a third-order tensor of size ${M_1^2} \times s \times 3$, where each slice is a covariance matrix. Calculate n smallest eigenvalue ${\lambda _1},{\lambda _{2}}, \cdots ,{\lambda _{n}}$ of the block diagonal matrix $bdiag(B)$ of this tensor.  Combining the model parameters $\theta _\sigma ^0,\theta _\sigma ^1, \cdot  \cdot  \cdot, \theta _\sigma ^n$ obtained from the training set, the final noise level estimation result is $\hat \sigma^2  = \theta _{{\sigma _0}}^0 + \theta _{{\sigma _0}}^1{\lambda _1} +  \cdots  + \theta _{{\sigma _0}}^n{\lambda _{n}}$.

\begin{algorithm}[h!]
\caption{Color image noise level estimation based on tensor decomposition}\label{algo1}
\hspace*{0.02in} {\bf Input:} Noisy image  ${\cal A}$;  Training set $ \cal \bar A$; termination distance $\varepsilon $\\
\hspace*{0.02in} {\bf Output:} Noise level  $\hat \sigma$ 
\begin{algorithmic}[1]
\State  ${W_{(k)}}$ $\Leftarrow$ $M_1 \times M_2 \times 3$  overlapping patches of  $ \cal \bar A$
\State Compute   $bdiag(B)$ with the homogeneous block tensors in ${W_{(k)}}$
\State $\lambda _1^{(i)},\lambda _{2}^{(i)}, \cdot  \cdot  \cdot ,\lambda _{n-1}^{(i)},\lambda _{n}^{(i)}(i = 1,2, \cdot  \cdot  \cdot ,M)$ $\Leftarrow$  eigenvalue of $bdiag(B)$  
\State  $train\_x$  $\Leftarrow$  $\lambda _1^{(i)},\lambda _{2}^{(i)}, \cdot  \cdot  \cdot ,\lambda _{n-1}^{(i)},\lambda _{n}^{(i)}(i = 1,2, \cdot  \cdot  \cdot ,M)$
\State $train\_y$ $\Leftarrow$ $\sigma  \pm \varepsilon $
\State $\theta _\sigma ^0,\theta _\sigma ^1, \cdot  \cdot  \cdot, \theta _\sigma ^n$  $\Leftarrow$ learn from $train_x$ and $train_y$ 
\State $w_{(k)}$ $\Leftarrow$ $M_1 \times M_2 \times 3$  overlapping patches of  ${\cal A}$
\State Compute   $bdiag(B)$ with the homogeneous block tensors in ${w_{(k)}}$
\State ${\lambda _1},{\lambda _{2}}, \cdots ,{\lambda _{n}}$ and ${\sigma _0}$ $\Leftarrow$ eigenvalue of $bdiag(B)$  
\State $\hat \sigma^2  = \theta _{{\sigma _0}}^0 + \theta _{{\sigma _0}}^1{\lambda _1} +  \cdots  + \theta _{{\sigma _0}}^n{\lambda _{n}}$ 
%
\end{algorithmic}
\end{algorithm}

\section{Experimental result }\label{4}

In this section, numerical experiments were conducted in the color image database TID2008 $TID2008$ \cite{ref8}. And  compare with traditional noise level estimation algorithms. The experiments for determining the optimal parameters and model parameters are conducted on the training set in the database $BSD500$ \cite{ref10}. In the comparative experiment, all algorithms are implemented in the  environment of $Matlab{\text{ }}R2020a{\text{ }}{\text{ }}(Intel(R){\text{ }}UHD {\text{ }} Graphics)$.


\subsection{Accuracy measure}

The estimated value obtained by our noise level estimation algorithm is denoted as $\sigma_i$. We use the root mean square error($RMSE$) and the mean absolute error($MAE$) to evaluate the performance of the algorithm. the $RMSE$ can be decomposed into the following form:
\begin{equation}\label{M28}
RMSE = \sqrt {\frac{{\sum\nolimits_{i = 1}^n {{{({\sigma _i} - {{\bar \sigma }_i})}^2}} }}{N}},
\end{equation}
where ${{\bar \sigma }_i}$ represents the mean value of the estimations$ \sigma_i$, $N$ represents the total number of experimental samples. The smaller the value of $RMSE$ is, the higher the accuracy of the algorithm will be. The calculation formula for $MAE$ is as follows:
\begin{equation}\label{M29}
MAE = \frac{{\sum\nolimits_{i = 1}^n {\left| {{\sigma _i} - {{\bar \sigma }_i}} \right|} }}{N}.
\end{equation}
The smaller the value of MAE is, the smaller the error of the algorithm will be.

\subsection{Parameter determination }

\begin{figure}[htpb]
\centering
\vspace{-0.2cm}
\subfigure[]{
\includegraphics[width=80mm,height=50mm]{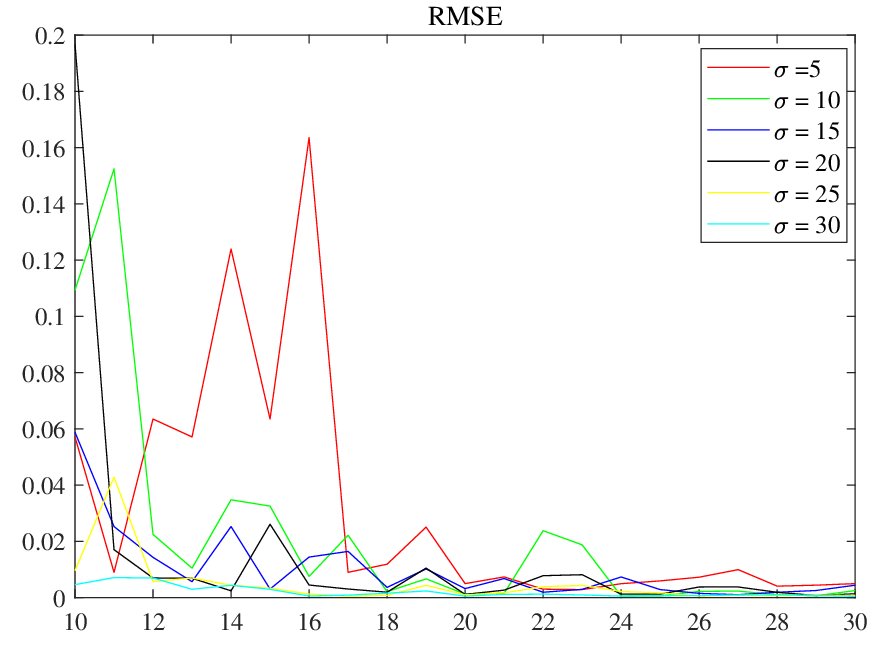}
}
\vspace{-2mm}
\subfigure[]{
\includegraphics[width=80mm,height=50mm]{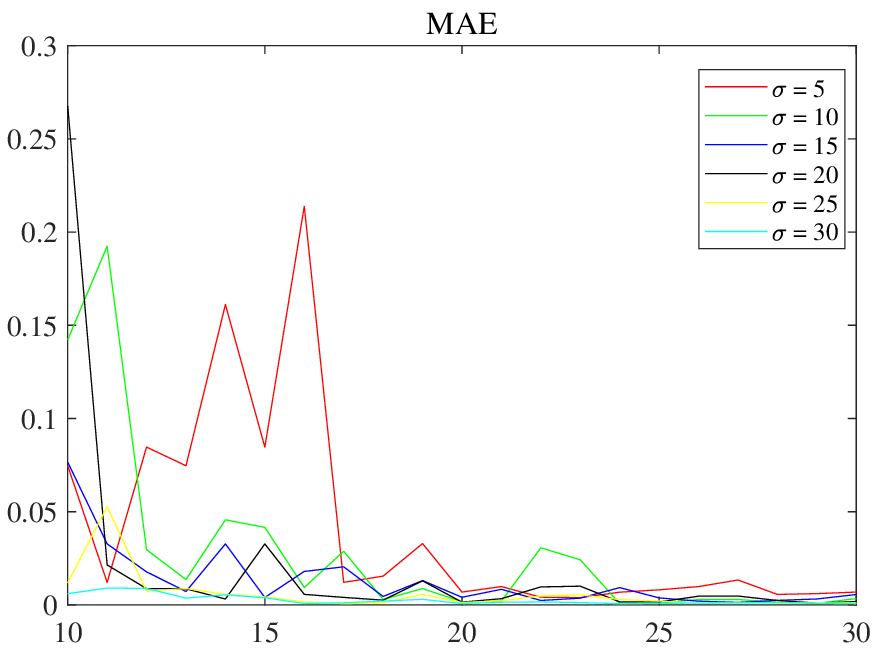}
}
\vspace{-2mm}

\caption{$RMSE$ and $MAE$ values corresponding to different $M$.}
\label{fig:1}
\end{figure}

There are three main parameters in the learning algorithm, namely$M_1$, $M$ and $n$.  $M_1$ represents the size of the selected sliding block. $M$ represents the number of images in the training set. $n$ represents the number of selected eigenvalues. In order to choose the optimal value, this paper use $RMSE$ as decision criteria. Parameter $M_1=7$ retains the value adopted by $Liu$ et al\cite{ref8}. Randomly selected $M$ training images in the $BSD500$ database  as the training set, the parameters $M$ of Eq. (\ref{M24}) experiment takes 10 to 30 in steps of one, respectively. The noise levels of the 200 images in $BSD500$ database with noise levels of 5, 10, 15, 20, 25, and 30 respectively are estimated. The estimated results are analyzed as shown in Fig. \ref{fig:1}. For overall performance, the $RMSE$ values that correspond to $M = 10$, 11, 12, 13, 14, 15, and 16 are larger than the $RMSE$ values that correspond to $M \geqslant 20$ in most cases. It is clear that the minimum estimation error obtained by the proposed algorithm with $M=20$, 21, 24, 29, and 30. However, the more images in training set, the longer the algorithm in this paper takes, so $M=20$ is selected in the following experiment to reduce the calculation time. 

\begin{figure}[htpb]
\centering
\vspace{-0.2cm}
\subfigure[]{
\includegraphics[width=80mm,height=50mm]{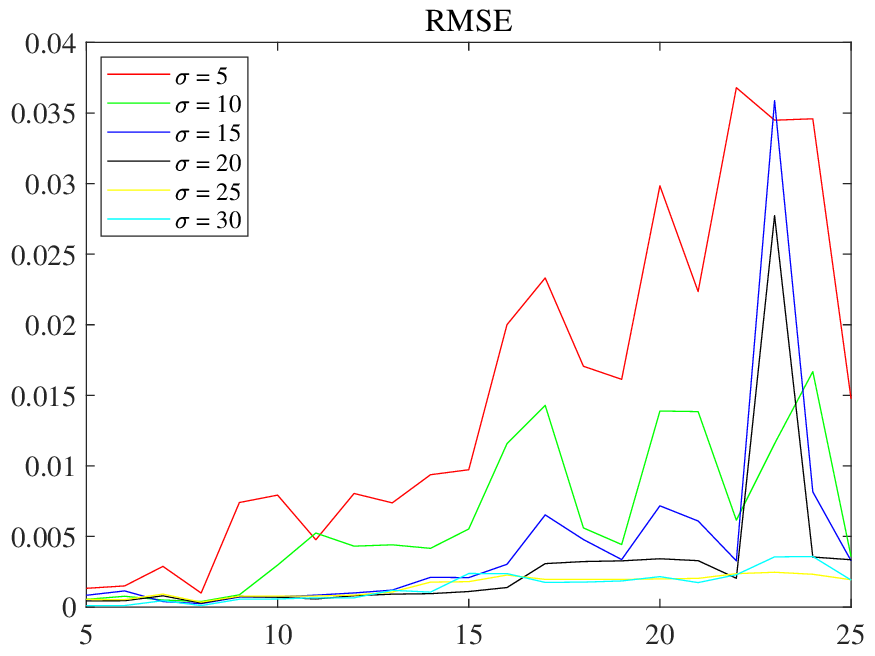}
}
\vspace{-2mm}
\subfigure[]{
\includegraphics[width=80mm,height=50mm]{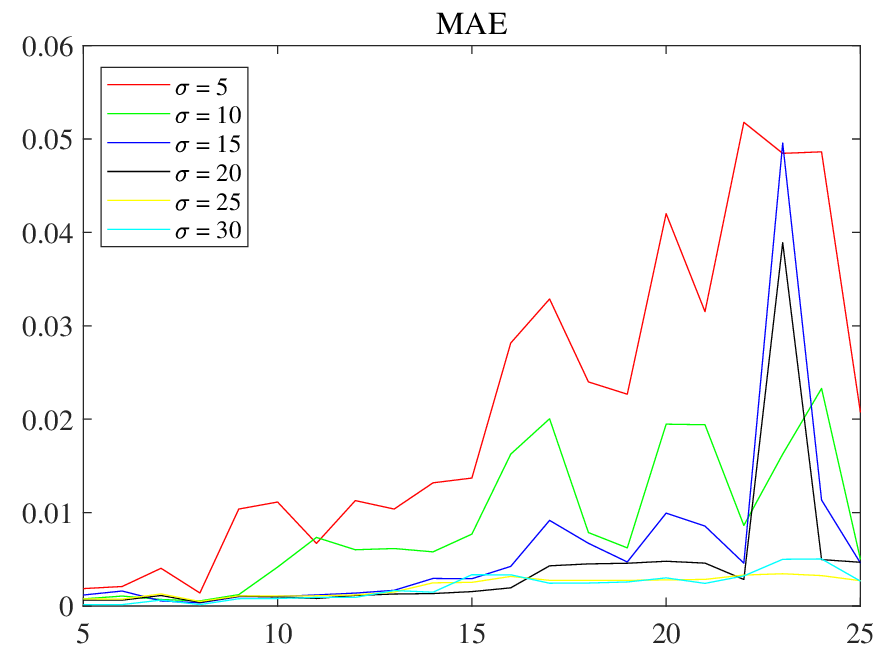}
}
\vspace{-2mm}

\caption{$RMSE$ and $MAE$ values corresponding to different $n$.}
\label{fig:2}
\end{figure}

The estimation effect of the algorithm is related to the selection of the number of eigenvalues $n$. This papertested parameter $n$ on 200 test set images in the database $BSD500$, and made the number of eigenvalues $n$ from 5 to 25 in steps of one to measure images processed by additive white Gaussian noise levels of 5, 10, 15, 20, 25, and 30.   The $MAE$ and $Bias$ values were calculated with the obtained results respectively, as shown in Fig. \ref{fig:2}. 

Abscissa for the value of $n$ in Fig. \ref{fig:2}, as the number of selected eigenvalues increases, the values of $RMSE$ and $MAE$ also increase. When the value of $n$ is relatively small, we can observe that at different levels of the noise levels, when $n=8$, the values of $RMSE$ and $MAE$ are both relatively small. So $n=8$ is selected in the following experiments.

\subsection{The estimation effect of the proposed algorithm }

We randomly selected 20 images with known noise levels from database $BSD500$ as the training set, and the learning parameters obtained are shown in Table \ref{tab1}. From Table 1, there are different learning coefficients under different noise levels, and this algorithm has relatively good estimation results.

\begin{table}[htbp]
   
\footnotesize
\begin{center}
  \centering
  \caption{Learning coefficients and estimation result.}
    \label{tab1}%
    \begin{tabular}{cccccccccc}
    \toprule
    \multicolumn{1}{c}{${\theta _0}$} & \multicolumn{1}{c}{${\theta _1}$} & \multicolumn{1}{c}{${\theta _2}$} & \multicolumn{1}{c}{${\theta _3}$} & \multicolumn{1}{c}{${\theta _4}$} & \multicolumn{1}{c}{${\theta _5}$} & \multicolumn{1}{c}{${\theta _6}$} & \multicolumn{1}{c}{${\theta _7}$} & \multicolumn{1}{c}{${\theta _8}$}& \multicolumn{1}{c}{result} \\
    \midrule
24.7839 & 0.0436 & -0.0674 & -0.1561 & 0.2159 & -0.0686 & 0.0079 & -0.1645 & 0.0646 & 5.0015\\
100.7873 & -0.0343 & 0.0154 & -0.0408 & -0.0355 & 0.0031 & -0.0453 & 0.0436 & 0.0039 & 9.9994\\
225.6460 & -0.0189 & 0.0142 & 0.0130 & -0.0088 & 0.0009 & -0.0091 & 0.0050 & 0.0009 & 14.9994\\
399.8086 & 0.0065 & -0.0038 & -0.0069 & -0.0039 & 0.0094 & -0.0026 & 0.0061 & -0.0043 & 19.9996\\
625.1615 & -0.0004 & 0.0043 & 0.0060 & -0.0087 & -0.0028 & 0.0111 & -0.0159 & 0.0061 & 25.0003\\
900.1740 & -0.0013 & 0.0026 & -0.0094 & 0.0099 & 0.0021 & -0.0042 & 0.0022 & -0.0019 & 30.0011\\
    \bottomrule
    \end{tabular}%

    \end{center}

\end{table}%

To further verify the feasibility of this algorithm.  Experiments have been carried out in the image database $TID2008$, and the results are shown in Table  \ref{tab2}. Under various noise levels, the algorithm presented in this paper can always obtain accurate estimation values. As can be seen from the table \ref{tab2}, the maximum error of the estimated result obtained by this algorithm is 0.0075 for the $'8.bmp'$ when the noise level is 20. This result has been presented in bold and italic. The minimum error can be as low as 0.0001.

\begin{table}[htbp]
\footnotesize
  \centering
  \caption{The result of noise estimation is obtained by the 25 images in image database $TID2008$ with additive white Gaussian noise of 5, 10, 15, 20, 25, and 30 respectively.}
    \begin{tabular}{lcccccc}
    \hline
    \diagbox{image}{ noiselevel} & 5     & 10    & 15    & 20    & 25    & 30 \\
    \hline
    \centering
    1.bmp & 5.0014 & 10.0011 & 15.0001 & 20.0004 & 24.9981 & 30.0005 \\
    2.bmp & 5.0019 & 9.9984 & 15.0002 & 19.9997 & 25.0003 & 29.9998 \\
    3.bmp & 4.9991 & 9.9982 & 15.0013 & 19.9994 & 24.9996 & 29.9999 \\
    4.bmp & 5.0011 & 10.0001 & 15.0001 & 19.9997 & 24.9999 & 30.0001 \\
    5.bmp & 5.0004 & 9.9997 & 15.0039 & 20.0002 & 24.9995 & 29.9997 \\
    6.bmp & 5.0005 & 10.0021 & 15.0022 & 19.9998 & 24.9996 & 30.0004 \\
    7.bmp & 4.9997 & 9.9989 & 14.9999 & 19.9993 & 25.0003 &30.0002 \\
    8.bmp & 5.0007 & 9.9977 & 14.9994 & \textbf{20.0075} & 24.9991 & 29.9998 \\
    9.bmp & 5.0046 & 9.9986 & 14.9983 & 19.9998 & 25.0003 & 30.0001 \\
    10.bmp & 5.0005 & 9.9990 & 15.0003 & 19.9994 & 24.9995 & 30.0011 \\
    11.bmp & 5.0049 & 9.9996 & 14.9997 & 19.9992 & 25.0010 & 30.0003 \\
    12.bmp & 5.0208 & 10.0011 & 14.9998 & 19.9997 & 24.9978 & 30.0008 \\
    13.bmp & 5.0011 & 9.9998 & 15.0002 & 20.0001 & 25.0004 & 30.0017 \\
    14.bmp & 5.0052 & 9.9994 & 14.9997 & 20.0010 & 24.9997 & 30.0001 \\
    15.bmp & 4.9995 & 9.9982 & 14.9999 & 19.9997 & 25.0005 & 29.9996 \\
    16.bmp & 5.0057 & 9.9983 & 14.9990& 19.9986 & 24.9994 & 29.9916 \\
    17.bmp & 5.0011 & 9.9985 & 14.9994 & 19.9981 & 25.0012 & 29.9969 \\
    18.bmp & 5.0057 & 9.9981 & 15.0002 & 19.9994 & 25.0009 & 30.0006 \\
    19.bmp & 5.0008 & 9.9974 & 15.0003 & 19.9995 & 24.9998 & 29.9985 \\
    20.bmp & 5.0017 & 9.9992 & 14.9995 & 19.9999 & 25.0013 & 29.9982 \\
    21.bmp & 5.0019 & 9.9985 & 14.9995 & 20.0001 & 24.9987 & 29.9972 \\
    22.bmp & 5.0006 & 10.0005& 14.9989 & 19.9999 & 24.9991 & 29.9969 \\
    23.bmp & 4.9994 & 9.9975 & 15.0012 & 20.0003 & 24.9999 & 30.0011 \\
    24.bmp & 5.0032 & 10.0051 & 14.9994 & 19.9996 & 25.0003 & 30.0005 \\
  25.bmp & 5.0005 & 10.0001 & 14.9998 & 19.9998 & 25.0001 & 30.0002 \\
    \bottomrule
    \end{tabular}%
  \label{tab2}%
\end{table}%

\subsection{Experimental comparison  }
The algorithm in this paper was compared with the traditional algorithms such as $Pyatykh$\cite{ref7} and $Liu$\cite{ref8}, as well as the latest algorithms such as $Fang$\cite{ref10} and $Liu$\cite{ref12}. Table \ref{tab3} presents the estimation results of 25 images in database $TID2008$ under different algorithms and at various noise levels. The optimal results have been marked in bold. From Table \ref{tab3}, at different noise levels, the estimation results of the algorithm in this paper were not the best in a few individual image experiments, but were the best in the results of the remaining image experiments.

\begin{table}[htbp]
\footnotesize

 \caption{The results of noise level estimation of color images with noise levels of 5, 10, 15, 20, 25, and 30 in database $TID2008$ are obtained by the following five algorithms.}
 (a)
  \centering
  
    \begin{tabular}{ccccccccccc}
    \toprule
     \diagbox{image}{ algorithms} &   $Pyatykh$    &   $Liu$    &   $Fang$   &   $Liu$    &   $Our$    &  $Pyatykh$    &   $Liu$    &   $Fang$     &   $Liu$    & $Our$  \\
\cmidrule{2-11}    \multicolumn{1}{r}{} & \multicolumn{5}{c}{\multirow{1}[2]{*}{    $\sigma  = 5$   }} & \multicolumn{5}{c}{\multirow{1}[2]{*}{$\sigma  = 10$}} \\
    \midrule
    1.bmp & 5.1777 & 4.9909 & 5.2812 & 5.0104 & \textbf{5.0014}  & 10.0276 & 10.0905 & 10.0653  & 10.0396 & \textbf{10.0011} \\
    2.bmp & 5.1238 & 4.9978 & 5.2790 & 4.9832 & \textbf{5.0019} & 10.0896 & 10.0095 & 10.1354 &10.0043 &  \textbf{9.9984} \\
    3.bmp & 5.1338 & 5.0321 & 5.0848 & 4.9497 & \textbf{4.9991} & 10.0011 & 10.0523 & 10.0221 & 10.0103 & \textbf{9.9982} \\
    4.bmp & 5.0920 & 5.1244 & 5.1326 & 4.9684 & \textbf{5.0011} & 10.1930 & 10.0084 & 10.1869 & 10.0067 & \textbf{10.0001} \\
    5.bmp & 5.1531 & 5.0496 & 5.2847 & 5.0374 & \textbf{5.0004} & 9.9989 & 9.9237 & 9.9977 & 10.0035 & \textbf{9.9997} \\
    6.bmp & 5.0561 & 4.9822 & 5.1409  & 5.0676 & \textbf{5.0005} & 10.1566 & 10.0169 & 10.1011 & 10.0201 & \textbf{10.0021} \\
    7.bmp & 5.1688 & \textbf{5.0002} & 5.1893 & 4.9787 & 4.9997 & 10.1732 & 9.9878 & 10.0268 & 9.9921 & \textbf{9.9989} \\
    8.bmp & 5.1261 & 5.2922 & 5.1158 & 5.3214 & \textbf{5.0007} & 10.0354 & 10.1505 & 10.3367 &10.0394 &  \textbf{9.9977} \\
    9.bmp & 5.2251 & 5.1719 & 5.1487 & 4.9737 & \textbf{5.0046} & 10.0370 & 10.0240 & 10.1496 & \textbf{9.9993} & 9.9986 \\
    10.bmp & 5.3007 & 5.1414 & 5.2805 & 5.0457 & \textbf{5.0005} & 10.0059 & 10.0194 & 10.0194 & 10.0042 & \textbf{9.9990} \\
    11.bmp & 5.1112 & 5.0503 & 5.2445 & 5.0124 & \textbf{5.0049} & 10.0023 & 9.9884 & 10.0177 & 10.0019 & \textbf{9.9996} \\
    12.bmp & 5.1375 & 5.1116 & 5.1377 & 5.0340 & \textbf{5.0208} & 10.0976 & 10.0119 & 10.1131 & \textbf{9.9990} & 10.0011 \\
    13.bmp & 5.2025 & 5.0151 & 5.2595 & 5.2659 & \textbf{5.0011} & 10.0874 & 10.2905 & 10.1191 & 10.0038 & \textbf{9.9998} \\
    14.bmp & 5.1953 & 5.0058 & 5.0231 & 5.1227 & \textbf{5.0052} & 10.057 & 9.9833 & 10.1429 & 9.9974 & \textbf{9.9994} \\
    15.bmp & 5.2204 & 5.0644 & 5.1572 & 5.0298 & \textbf{4.9995} & 10.1362 & 10.0295 & 10.1689 & \textbf{10.0016} & 9.9982 \\
    16.bmp & 5.1995 & 5.0320 & 5.0539 & 4.9429 & \textbf{5.0057} & 10.0845 & 10.0256 & 10.0845 & 10.0048 & \textbf{9.9983} \\
    17.bmp & 5.2315 & 5.1049 & 4.9607 & 5.0593 & \textbf{5.0011} & 10.3109 & 9.9945 & 10.1253 & 10.0128 & \textbf{9.9985} \\
    18.bmp & 5.0488 & 5.2469 & 5.0784  & 5.0822 & \textbf{4.9985} & 10.1735 & 9.9962 & 10.1483 & 9.9956 & \textbf{9.9981} \\
    19.bmp & 5.1256 & 5.0760 & 5.0733 & 4.9620 & \textbf{5.0008} & 10.0938 & 10.0048 & 10.1653 & \textbf{10.0026} & 9.9973 \\
    20.bmp & 5.0614 & 4.9974 &5.0784 & 4.9688 & \textbf{5.0017} & 9.9381 & 9.9791 & 10.0780 & \textbf{9.9997} & 9.9992\\
    21.bmp & 5.2948 & 5.1776 & 4.9888 & 5.0884 & \textbf{5.0019} & 10.2443 & 9.9552 & 10.1925 & 10.0075 & \textbf{9.9985} \\
    22.bmp & 5.2739 & 5.1267 & 5.1268 & 5.0235 & \textbf{5.0006} & 10.1127 & 10.0087 & 10.1338 & 9.9938 & \textbf{10.0005} \\
    23.bmp & 5.2212 & 5.0940 & 5.0672 & 4.9614 & \textbf{4.9994} & 10.0801 & 10.0010 & 10.1001 & \textbf{9.9999} & 9.9975 \\
    24.bmp & 5.0355 & 4.9980 & 5.1597 & 5.0269 & \textbf{5.0026} & 10.1337 & 9.9974 & 10.0570 & 10.0059 & \textbf{10.0051} \\
    25.bmp & 5.0565 & 4.9061 & 4.9749 & 4.9289 & \textbf{5.0005} & 10.1134 & 9.8087 & 9.9729 & 9.9982 & \textbf{10.0001} \\
    \bottomrule
    \end{tabular}%
      \label{tab3}%
\end{table}%
  \begin{table}[htbp]
\footnotesize
(b)
\centering

  \begin{tabular}{ccccccccccc}
    \toprule
     \diagbox{image}{ algorithms} &   $Pyatykh$    &   $Liu$    &   $Fang$   &   $Liu$    &   $Our$    &  $Pyatykh$    &   $Liu$    &   $Fang$     &   $Liu$    & $Our$  \\
\cmidrule{2-11}    \multicolumn{1}{r}{} & \multicolumn{4}{c}{\multirow{1}[2]{*}{    $\sigma  = 15$   }} & \multicolumn{4}{c}{\multirow{1}[2]{*}{$\sigma  = 20$}} \\
    \midrule
    1.bmp & 14.8888 & 14.9418 & 15.2009 & 15.0761 & \textbf{15.0001} & 20.4419 & 20.0425 & 20.3151 & 19.9899 & \textbf{20.0004} \\
    2.bmp & 15.1404 & 14.9667 & 15.0645 & 15.0062 & \textbf{15.0002} & 19.9699 & 19.9743 & 19.9991 & 20.0016 & \textbf{19.9997} \\
    3.bmp & 15.0828 & 14.8765 & 14.9452 & 14.9952 & \textbf{15.0013} & 19.8293 & 19.8108 & 19.9986 & 20.0011 & \textbf{19.9994} \\
    4.bmp & 14.9683 & 14.9402 & 15.1621 & 14.9936 & \textbf{15.0001} & 19.9252 & 19.9568 & 19.9989 &  19.9993 & \textbf{19.9997} \\
    5.bmp & 14.9465 & 14.9728 & 15.0824 & 14.9918 & \textbf{15.0039} & 20.1196 & 19.8898 & 19.9909 & 20.0017 & \textbf{20.0002} \\
    6.bmp & 15.0655 & 15.0085 & \textbf{14.9999} & 15.0099 & 15.0022 & 20.1091 & 19.9724 & 20.1928 & 19.9974 & \textbf{19.9998} \\
    7.bmp & 15.0667 & 15.0006 & 15.0665 & 15.0039 & \textbf{14.9999} & 20.0524 & 19.8218 & 20.0119 & 19.9946 & \textbf{19.9993} \\
    8.bmp & 15.0554 & 15.0088 & 15.2302 & 14.9904 & \textbf{14.9994} & 20.1131 & 20.0678 & 20.1523 & \textbf{20.0012} & 20.0075 \\
    9.bmp & 14.9893 & 14.9909 & 15.0438 & 15.0049 & \textbf{14.9983} & 19.9822 & 19.9696 & 20.0541 & \textbf{20.0009} & 19.9983 \\
    10.bmp & 15.0197 & 14.9917 & 15.0868 & 15.0229 & \textbf{15.0003} & 20.0536 & 19.8758 & 20.0953 & 19.9959 & \textbf{19.9994} \\
    11.bmp & 15.0039 & 14.9957 & 15.0875 & 14.9884 & \textbf{14.9997} & 20.0842 & 20.0134& 20.1256  & 20.0018 & \textbf{19.9992} \\
    12.bmp & 15.0151 & 15.0487 & 15.1621 & 15.0017 & \textbf{14.9998} & 19.9609 & 19.8579 & 20.0011 & 20.0013 & \textbf{19.9997} \\
    13.bmp & 14.9637 & 15.2757 & 15.1716 & 15.0023 & \textbf{15.0002} & 20.1668 & 20.4568 & 20.5263 & 20.0062 & \textbf{20.0001} \\
    14.bmp & 15.0943 & 15.0096 & 15.0730 & 15.0187 & \textbf{14.9997} & 20.1985 & 19.9788  & 20.1498 & 19.9927 & \textbf{20.0010} \\
    15.bmp & 15.1630 & 14.9740 & 15.1525 & 15.0122 & \textbf{14.9999} & 19.9105 & 20.0109 & 20.1471 & 19.9927 & \textbf{19.9997} \\
    16.bmp & 15.0667 & 14.9916 & 15.0854 & 14.9944 & \textbf{14.9990} & 20.0507 & 19.9725 & 20.0407 & 20.0029 & \textbf{19.9986} \\
    17.bmp & 15.1603 & \textbf{15.0004} & 15.0157 & 15.0033 & 14.9994 & 19.8769 & 19.9923 & 20.0504 & 19.9972 & \textbf{19.9981} \\
    18.bmp & 15.0636 & 14.9904  & 15.1878& 15.0263 & \textbf{15.0002} & 20.1855 & 20.0019 & 20.0929 & 19.9991 & \textbf{19.9994} \\
    19.bmp & 15.2471 & 14.9307 & 15.0729 & 14.9992 & \textbf{15.0003} & 20.2273 & 19.9720 & 20.2322 & 19.9955 & \textbf{19.9995} \\
    20.bmp & 14.6112  & 14.9864 & 15.1039 & 14.9967 & \textbf{14.9995} & 19.6693 & 19.9619 & 20.0987 & 20.0052 & \textbf{19.9999} \\
    21.bmp & 15.1656 & 14.9594 & 15.0728 & 15.0152 & \textbf{14.9995} & 20.1306 & 19.9663 & 19.9749 & 19.9945 & \textbf{20.0001} \\
    22.bmp & 15.2042 & 14.9082 & 15.1112 & 14.9988 & \textbf{14.9989} & 20.4745 & 20.0349 & 20.1887 & 19.9998 & \textbf{19.9999} \\
    23.bmp & 14.9134 & 14.9316 & 14.9848 & 14.9966 & \textbf{15.0012} & 19.7002 & 19.9387 & 20.1008 & 20.0050 & \textbf{20.0003} \\
    24.bmp & 15.3326 & 14.9815 & 15.1788 & \textbf{15.0005} & 14.9994 & 20.2958 & 20.0202 & 19.9172 & 20.0011 & \textbf{19.9996} \\
    25.bmp & 14.8590 & 14.8689 & 14.9701 & 15.0016 & \textbf{14.9998} & 19.9431 & 19.7954 & 19.8702  & 20.0072 & \textbf{19.9998} \\
    \bottomrule
    \end{tabular}%
\end{table}%

\begin{table}[htbp]
\footnotesize
(c)
\centering

  \begin{tabular}{ccccccccccc}
   \toprule
     \diagbox{image}{ algorithms} &   $Pyatykh$    &   $Liu$    &   $Fang$   &   $Liu$    &   $Our$    &  $Pyatykh$    &   $Liu$    &   $Fang$     &   $Liu$    & $Our$  \\
\cmidrule{2-11}    \multicolumn{1}{r}{} & \multicolumn{5}{c}{\multirow{1}[2]{*}{    $\sigma  = 25$   }} & \multicolumn{5}{c}{\multirow{1}[2]{*}{$\sigma  = 30$}} \\
    \midrule
        1.bmp & 25.1796 & 25.0353 & 25.2086 & 25.0112 & \textbf{24.9981} & 30.1183 & 30.0873 & 30.2582 & 29.9969 & \textbf{30.0005} \\
    2.bmp & 24.9384 & 24.8494 & 25.2505 & 24.9949 & \textbf{25.0003} & 29.9177 & 29.8028 & 30.0896 & 29.9996 & \textbf{30.0003} \\
    3.bmp & 24.6657 & 24.8306 & 25.0050 & 24.9873 & \textbf{24.9996} & 28.9310 & 29.9016 & 30.0564 & 30.0048 & \textbf{29.9988} \\
    4.bmp & 24.7732 & 24.8843 & 25.0020 & 24.9915 & \textbf{24.9999} & 29.0958 & 29.7832 & 30.0556 & 29.9997 & \textbf{30.0001} \\
    5.bmp & 25.1836 & 25.0257 & 25.1081 & 25.0374 & \textbf{24.9995} & 30.3305 & 29.9751 & 30.1791 & \textbf{29.9998} & 29.9997 \\
    6.bmp & 25.1349 & 24.9826 & 25.1221 & 25.0156 & \textbf{24.9996} & 30.0484 & 29.9161 & 30.1003 & 29.9898 & \textbf{30.0004} \\
    7.bmp & 24.9195 & 24.8384 & 24.8583 & 24.9857 & \textbf{25.0003} & 30.0766 & 29.8047 & 29.7972 & 29.9981 & \textbf{30.0002} \\
    8.bmp & 24.8943 & 24.9657 & 25.1712 & 24.8761 & \textbf{24.9991} & 30.0125 & 29.9507 & 30.2289 & 29.9943 & \textbf{29.9998} \\
    9.bmp & 24.9576 & 24.8822 & 25.1047 & 24.9842 & \textbf{25.0003} & 29.9645 & 29.9439 & 30.0681 & 29.9991 & \textbf{30.0001} \\
    10.bmp & 24.9955 & 24.9320 & 24.9399 & 24.9728 & \textbf{24.9995} & 29.7311 & 29.8181 & 29.8880 & \textbf{29.9992} & 30.0011 \\
    11.bmp & 25.1671 & 25.0357 & 25.0850 & 24.9943 & \textbf{25.0010} & 30.0311 & 29.9151 & 30.2016 & 29.9949 & \textbf{30.0003} \\
    12.bmp & 25.0011 & 24.8305 & 25.0314 & 24.9718 & \textbf{24.9978} & 29.8704 & 29.6872 & 29.8480 & 29.9946 & \textbf{30.0008} \\
    13.bmp & 25.0522 & 25.3657 & 25.7569 & 25.0071 & \textbf{25.0004} & 30.0921 & 30.1268 & 30.2795 & 29.9901 & \textbf{30.0017} \\
    14.bmp & 25.0454 & 24.9507 & 25.0122 & 24.9981 & \textbf{24.9997} & 30.1545 & 29.9122 & 30.1025 & 29.9958 & \textbf{30.0001} \\
    15.bmp & 24.6330 & 24.9849 & 25.1236 & 24.9779 & \textbf{25.0005} & 29.9529 & 30.0174 & 30.0184 & 30.0029 & \textbf{29.9996} \\
    16.bmp & 24.9197 & 25.0613 & 24.9661 & 24.9807 & \textbf{24.9994} & 30.0878 & 29.9047 & 30.0754 & \textbf{29.9978} & 29.9916 \\
    17.bmp & 24.3706 & 24.9467 & 2\textbf{5.0009 }& 24.9873 & 25.0012 & 29.1384 & 29.7090 & 30.0498 & 29.9909 & \textbf{29.9969} \\
    18.bmp & 25.0396 & 24.9923 & 25.0849 & 24.9931 & \textbf{25.0009} & 30.2671 & 30.0018 & 30.1984 & 29.9972 & \textbf{30.0006} \\
    19.bmp & 24.9352 & 24.8534 & 25.2202 & 25.0314 & \textbf{24.9998} & 30.3306 & 29.9253 & 29.8960 & 29.9958 & \textbf{29.9985} \\
    20.bmp & 24.6186 & 24.8590 & 24.9886 & 24.9950 & \textbf{25.0012} & 29.7134 & 29.9794 & 29.9851 & 30.0021 & \textbf{29.9982} \\
    21.bmp & 25.0496 & 24.9992 & 25.0447 & \textbf{25.0012} & 24.9987 & 30.0134 & 30.0294 & 30.0197 & \textbf{30.0023} & 29.9972 \\
    22.bmp & 25.2642 & 25.0014 & 25.2647 & 24.9899 & \textbf{24.9991} & 30.0395 & 30.0434 & 30.1528 & 29.9957 & \textbf{29.9969} \\
    23.bmp & 24.3456 & 24.8050 & 24.9526 & 24.9773 & \textbf{24.9999} & 29.3856 & 29.7548 & 29.8333 & 29.9957 & \textbf{30.0011} \\
    24.bmp & 25.2909 & 24.7809 & 24.9411 & 25.0039 & \textbf{25.0003} & 30.1711 & 29.9615 & 30.1141 & 30.0214 & \textbf{30.0005} \\
    25.bmp & 25.2844 & 24.8446 & 24.9290 & 24.9518 & \textbf{25.0001} & 30.0902 & 30.0970 & 30.0543 & 29.9973 &  \textbf{30.0002} \\
    \bottomrule
    \end{tabular}%
\end{table}

\begin{figure}[htpb]
\centering
\vspace{-0.2cm}
\subfigure[]{
\includegraphics[width=45mm,height=40mm]{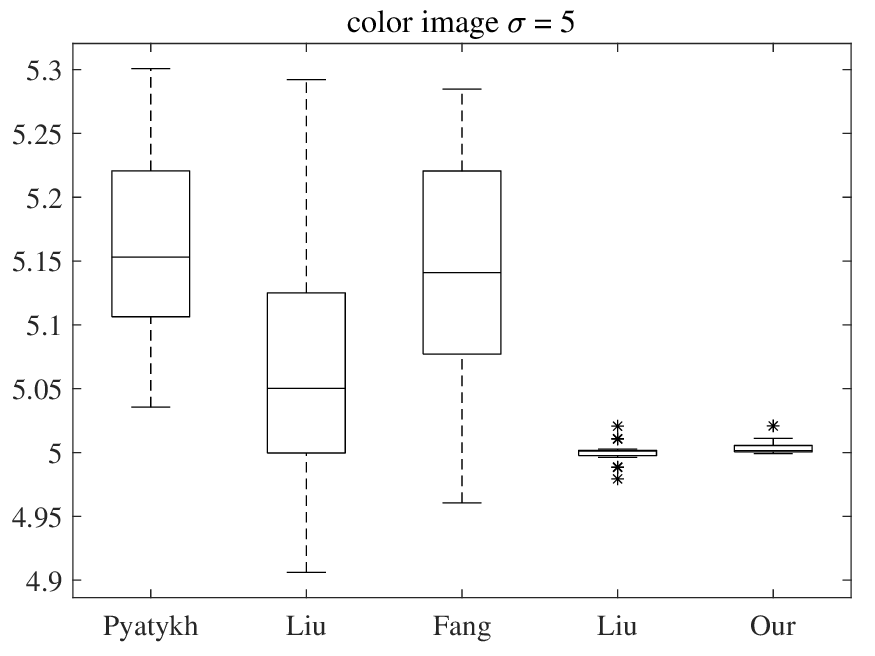}
}
\vspace{-2mm}
\subfigure[]{
\includegraphics[width=45mm,height=40mm]{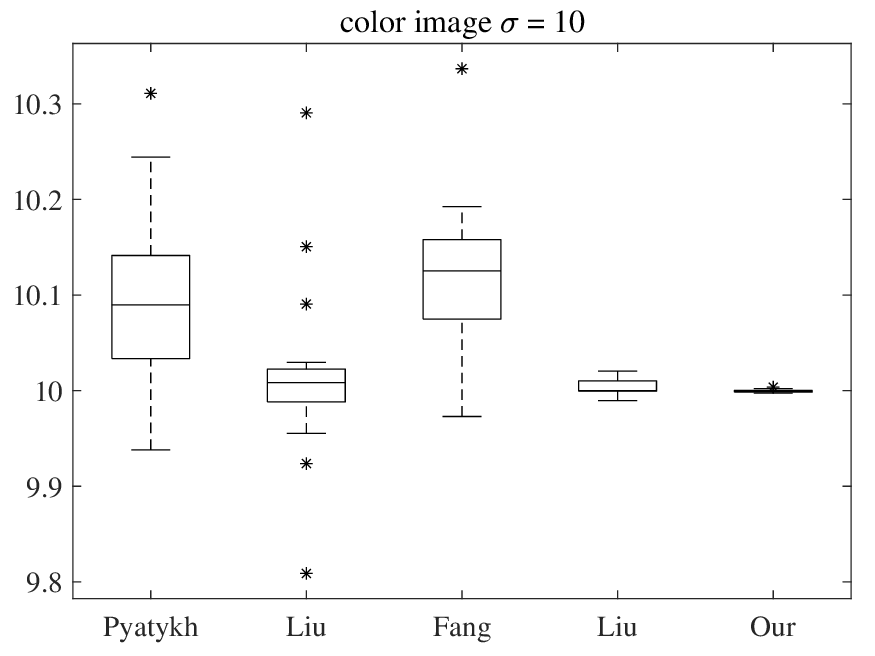}
}
\vspace{-2mm}
\subfigure[]{
\includegraphics[width=45mm,height=40mm]{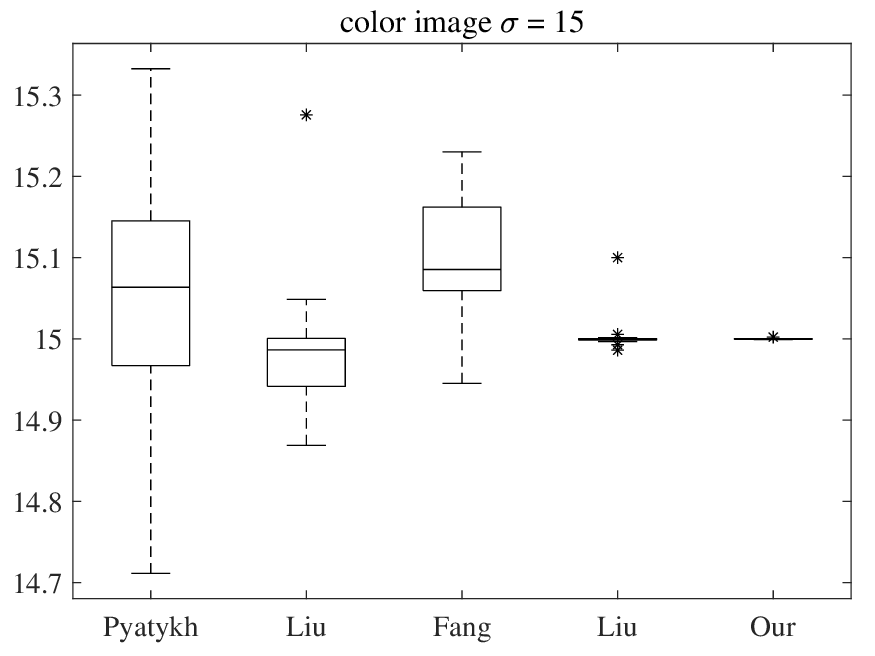}
}
\vspace{-2mm}
\subfigure[]{
\includegraphics[width=45mm,height=40mm]{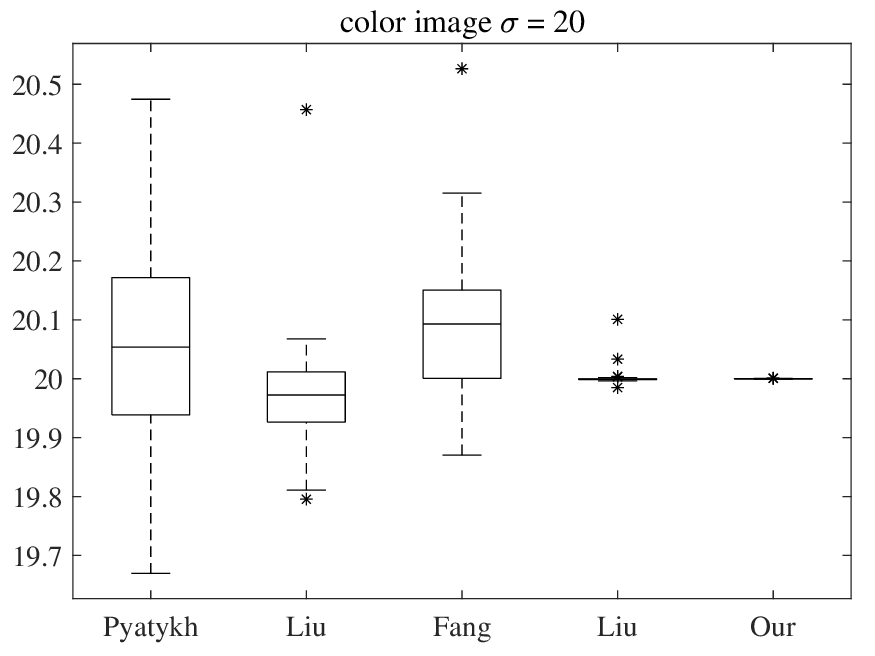}
}
\vspace{-2mm}
\subfigure[]{
\includegraphics[width=45mm,height=40mm]{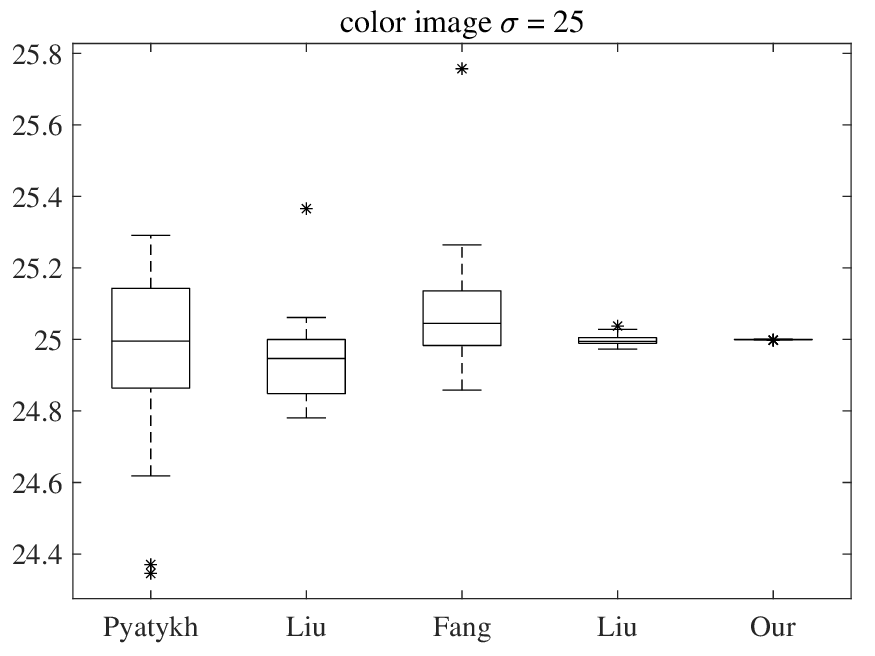}
}
\vspace{-2mm}
\subfigure[]{
\includegraphics[width=45mm,height=40mm]{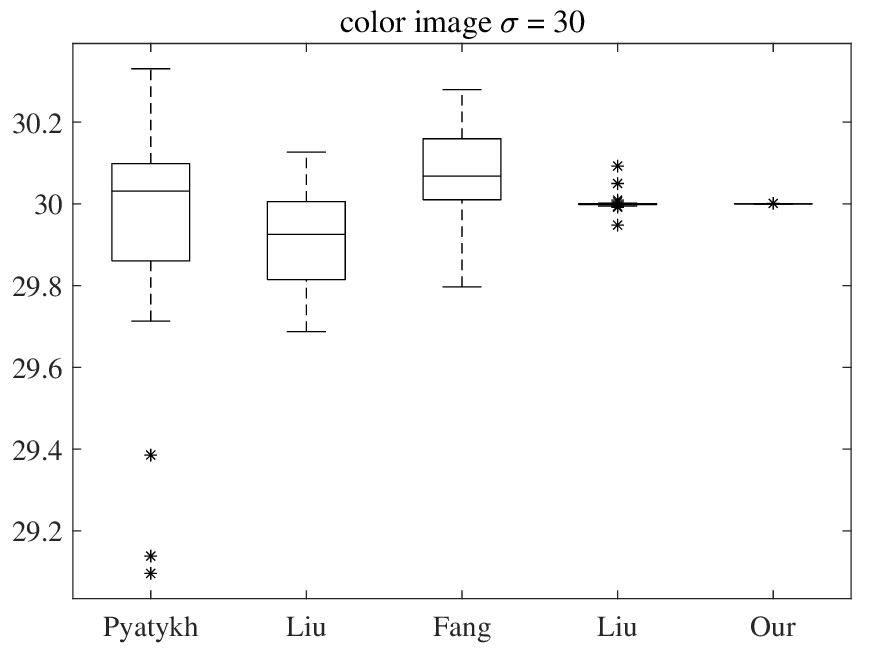}
}
\vspace{2mm}
\caption{The results of noise level estimation with noise levels of 5, 10, 15, 20, 25, and 30 are obtained by the five compared algorithms. }
\label{fig:3}
\end{figure}

To facilitate a more intuitive comparison of the estimation results of the five algorithms, we have drawn a box plot as shown in Fig. \ref{fig:3}. From Fig. \ref{fig:3}, because the methods in this paper and those of $Liu$ et al. were trained using the learning approach, the accuracy of these two methods is generally quite high. When the noise level is 15, 20 and 30, the accuracy of the method in this paper is similar to that of $Liu$ et al.'s method, and the number of outliers in the method of this paper is smaller. When the noise level is 5, 10 and 25, compared with the algorithm of $Liu$ et al., the fluctuation of the estimation results of the algorithm in this paper is smaller. This indicates that the algorithm in this paper has good stability.

\begin{figure}[htpb]
\centering
\vspace{-0.2cm}
\subfigure[]{
\includegraphics[width=80mm,height=50mm]{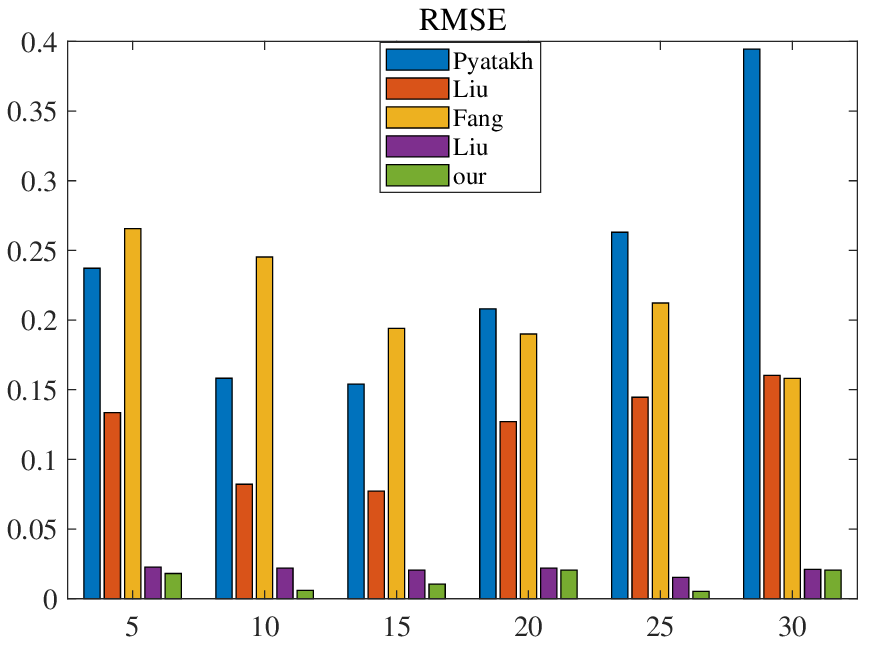}
}
\vspace{-2mm}
\subfigure[]{
\includegraphics[width=80mm,height=50mm]{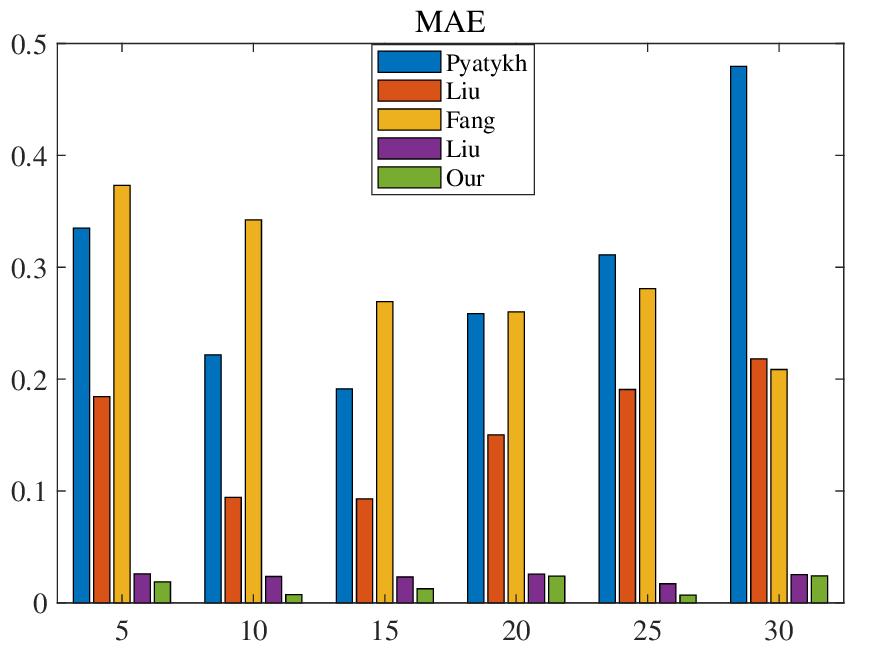}
}
\vspace{-2mm}

\caption{$RMSE$ and $MAE$ values corresponding to different algorithm.}
\label{fig:4}
\end{figure}

We present the results of $RMSE$ and $MAE$ as shown in Fig. \ref{fig:4}. The smaller the value of $RMSE$ and $MAE$ is, the higher the accuracy of the algorithm will be. When the noise level exceeds 20, the $RMSE$ and $MAE$ values of the $Pyatakh$ et al. algorithm reach their maximum. When the noise level is less than 20, the $RMSE$ and $MAE$ values of the $Fang$ et al. algorithm reach their maximum. When the noise level is 20 and 30, the $RMSE$ and $MAE$ of the algorithm in this paper are similar to those of the algorithm proposed by $Liu$ et al.. Overall, the algorithm presented in this paper has a relatively high estimation accuracy. 

Under the same operating conditions, we compared the running times of five algorithms, as shown in Table \ref{tab4}. From Table \ref{tab4}, the algorithm proposed by $Liu$ et al. has the shortest running time because they only conducted experiments using weakly textured graphics, which significantly improved the efficiency of the algorithm. Because the algorithm in this paper requires training the learning parameters in the training set, and decomposing the tensors corresponding to color images, the running time of this algorithm is relatively long.

\begin{table}[htbp]
\footnotesize
  \centering
  \caption{The result of noise estimation is obtained by the 25 images in image database $TID2008$ with additive white Gaussian noise of 5, 10, 15, 20, 25, and 30 respectively.}
    \begin{tabular}{lccccc}
    \hline
    \diagbox{time}{ algorithms} & $Pyatakh$     & $Liu$    & $Fang$    & $Liu$    & $Our$  \\
    \hline
    \centering
    $\min t(s)$ & 4.45 & 1.51 & 2.34 & 6.24 & 7.12  \\
    $\bar t(s)$ & 4.65 & 1.67 & 2.52 & 8.16 & 8.28  \\
    $\max t(s)$ & 4.77 & 2.16 & 2.67 & 9.64 & 9.82  \\
    \bottomrule
    \end{tabular}%
  \label{tab4}%
\end{table}%

\section{Conclusions and future work}\label{5}

This paper presents a method for estimating the noise level of color images based on tensor decomposition. The destruction of the color image tensor structure caused by traditional methods is avoided by this method, thereby improving the estimation accuracy of the algorithm. Through the proof of \textbf{Theorem 2}, it is demonstrated that the eigenvalues of the block diagonal matrix obtained by decomposing the tensor are related to the level of image noise. And through the learning method, the learning coefficient between the feature values and the noise level is trained. 

The parameters $M$(the number of images in the training set) and $N$(the number of eigenvalues) of the algorithm ware selected in database $BSD500$, and trained the learning parameters in the training set of $BSD500$. To illustrate the estimation accuracy of the algorithm in this paper, four algorithms ware compared in database $TID2008$. Under noise levels of 5, 10, 15, 20, 25 and 30, experiments were conducted on 25 color images in the database  $TID2008$. Plot the obtained results as a box plot, and calculate its $RMSE$ and $MAE$ values.  By comparison, it can be concluded that the algorithm presented in this paper has a relatively high estimation accuracy. Subsequently, we compared the running times of the five algorithms. Because the algorithm in this paper requires training of the learning parameters and decomposition and rearrangement of tensors, it is computationally intensive.

In this paper, this tensor decomposition method does not completely preserve the special structure of the tensor to a certain extent. The eigenvalues obtained through this method cannot fully represent the eigenvalues of the tensor. And the algorithm in this paper takes a relatively long time to execute. In future work, efforts will be made to maintain the structure of the tensor as much as possible, and the relationship between the corresponding eigenvalues of the tensor and the image noise level will be investigated. And shorten the running time of the algorithm.

\section*{Acknowledgments}This study was supported by  the National Natural Science Foundation of P.R. China (Grant No.12171064).


\section*{Declarations}

The authors declare no conflict of interest.
%
%

%
%

\end{document}